\documentclass[11pt,twoside]{amsart}
\usepackage{latexsym}
\usepackage{amssymb}
\newcommand{\eps}{\varepsilon}

\newcommand{\N}{{\mathbb N}}

\theoremstyle{plain}
\newtheorem{thm}{Theorem}[section]

\newtheorem{prop}[thm]{Proposition}
\newtheorem{cor}[thm]{Corollary}

\theoremstyle{definition}

\begin{document}
\title[Weak cluster points of a sequence and coverings by cylinders]%
{Weak cluster points of a sequence and coverings by cylinders}
\author{Vladimir Kadets}
\date{\today}
\address{Department of Mechanics and Mathematics,
Kharkov National University,\linebreak
pl.~Svobody~4, 61077~Kharkov, Ukraine}
\email{vova1kadets@yahoo.com}
\curraddr{Department of Mathematics, University of Missouri,
Columbia MO 65211}
\email{kadets@math.missouri.edu}
\thanks{I would like to say my thanks
to Nigel Kalton for his invitation to the University of Missouri,
where this research has been done, and
to Sveta Mayboroda, who asked the right question at the right time}
\begin{abstract}
Let $H$ be a Hilbert space. Using Ball's solution of the "complex
plank problem" we prove that the following properties of a sequence
$a_n>0$ are equivalent:
\begin{enumerate}
\item There is a sequence $x_n \in H$ with
$\|x_n\|=a_n$, having 0 as a weak cluster point;
\item
$\sum_1^\infty a_n^{-2}=\infty$.
\end{enumerate}
Using this result we show that a natural idea of generalization
of Ball's "complex plank" result to cylinders with $k$-dimensional
base fails already for $k=3$. We discuss also generalizations of
"weak cluster points" result
to other Banach spaces and relations with cotype.
\end{abstract}
\maketitle
\section{Introduction}
Let $H$ be an infinite-dimensional Hilbert space.
It is well known that the weak topology of $H$ has
bad sequential properties: a sequence $h_n \in H$
having a weak cluster point $x$, can be free from weakly
convergent to $x$ subsequences. Moreover, for a sequence $h_n \in H$
with a weak cluster point it is possible to have
$\|h_n\| \to \infty$ \cite{conkad}.
In this paper we study the following question:
if a sequence has a weak cluster point, how quick
can tend to infinity the norms of the sequence elements?
The main tool for estimation of this speed from above
will be Ball's "complex plank" theorem - a recent and
really beautiful statement from Hilbert space geometry.

Recall, that by a plank of width $w$ in $H$ one means a set of
the form
$$
P=\{h \in H : \left |\langle h - h_0 \, , \, e \rangle \right |
\le \frac{w}{2}\},
$$
where $\|e\| = 1$.
According to T.Bang's theorem \cite{bang}, if a sequence $P_n$ of planks
of widths $w_n$ covers a ball of diameter $w$ , then
$\sum w_n \ge w$.

Bang's theorem is true both in real and complex spaces, but in complex
spaces it can be dramatically improved. The following theorem is proved
by  K. Ball \cite{comp}
\begin{thm} \label{Ball}
Let
$$P_n=\{h \in H : \left |\langle h , e_n \rangle \right |
\le \frac{w_n}{2}\},$$
where $\|e_n\|=1$, be planks of widths $w_n$ in a complex
Hilbert space, and let
$\bigcup P_n $ cover a ball of diameter $w$ centered in origin.
Then $\sum w_n^2 \ge w^2$.
\end{thm}
We are going not only to apply the Ball's theorem to our problem,
but also to apply our result on weak topology to give some
limitations for possible generalizations of the Theorem \ref{Ball}.
After that we discuss generalizations
of our weak topology result to other Banach
spaces.
Allover the paper the letter $X$ will be used for infinite dimensional
Banach space, $S_X$ and $B_X$ - for its unit sphere and unit ball
respectively.
\section{Weighted weak convergence}

Let $P=\{p_{n,m}\}$ be an infinite matrix of non-negative numbers with
finite rows, satisfying conditions
\begin{enumerate}
\item
$\sum_{m=1}^\infty p_{n,m} = 1$.
\item $\lim_{n \to \infty} p_{n,m} =0$
\end{enumerate}
We say that a sequence $\{a_n\}$
of reals is $P$-convergent to 0 ($P\lim a_n =0$) if
$$
\lim_{n \to \infty} \sum_{m=1}^\infty p_{n,m}a_m = 0.
$$
Some evident well-known properties of $P$-convergence
are collected in the proposition below.
\begin{prop} \label{p-conv}
\begin{enumerate}
\item If \, $\lim a_n =0$, then $P\lim a_n =0$.
\item \label{p2}
If $P\lim a_n =0$ and $P\lim b_n =0$, then $P\lim (a_n + b_n)=0$. 
\item \label{p3}
If $a_n \ge 0$ and $P\lim a_n =0$, then 0 is a cluster point
of $\{a_n\}$.
\end{enumerate}
\end{prop}
Let $X$ be a Banach space, $1 \le p < \infty$.
A sequence $x_n \in X$ is said to be weakly $(P,p)$-convergent to 0
($w(P,p)\lim x_n =0$), if for every $f \in X^*$ the sequence
$ \left |f(x_n) \right |^p$ $P$-converges to 0.
\begin{prop} \label{wp-conv}
If $w(P,p)\lim x_n =0$, then 0 is a weak cluster point
of $\{x_n\}$.
\end{prop}
\begin{proof}
We must prove that for arbitrary finite set $\{f_k\}_{k=1}^m \subset X^*$ 
and arbitrary $\eps > 0$ there is an arbitrarily big $n \in \N$ with
$\sum_{k=1}^m \left | f_k (x_n) \right |^p \le \eps$.
By definition of $w(P,p)$-convergence and item (\ref{p2}) of
Proposition \ref{p-conv}
$$
P\lim_n \sum_{k=1}^m \left | f_k (x_n) \right |^p = 0.
$$
The rest we deduce from item (\ref{p3}) of
Proposition \ref{p-conv}.
\end{proof}

\section{The main result}
\begin{thm} \label{main}
The following properties of a sequence
$a_n>0$ are equivalent:
(a) \ \ \ $\sum_1^\infty a_n^{-2}=\infty$.
(b) \ There is a sequence $x_n$ in a Hilbert space $H$ with
$\|x_n\|=a_n$, and there is a sequence $P=\{p_{n,m}\}$ of weights
such that $w(P,2)\lim x_n =0$.
(c) \ There is a sequence $x_n$ in a Hilbert space $H$ with
$\|x_n\|=a_n$, having 0 as a weak cluster point.
\end{thm}
\begin{proof}
(a) $\Longrightarrow$ (b). Assume
$\sum_1^\infty a_n^{-2}=\infty$ and consider
the sequence $x_n=a_n e_n$, where $e_n$ is an orthonormal sequence
in $H$. Introduce the following sequence $P=\{p_{n,m}\}$ of weights:
$$p_{n,m}=\frac{a_m^{-2}}{\sum_{j=1}^n a_j^{-2}}$$
when $m \le n$ and $p_{n,m}=0$ for $m > n$.
Then for every $f \in H$ we have
$$
\sum_{m=1}^\infty p_{n,m} \left |\langle x_m, f \rangle \right |^2=
\frac{\sum_{m=1}^n
\left |\langle e_m, f \rangle \right |^2}{\sum_{m=1}^n a_m^{-2}} \le
\frac{\|f\|^2}{\sum_{m=1}^n a_m^{-2}} \to 0, \ as \ n \to \infty
$$
which means $w(P,2)$-convergence of $x_n$ to 0.
(b) $\Longrightarrow$ (c). This is given by the
Proposition \ref{wp-conv}.
(c) $\Longrightarrow$ (a). Let
\begin{equation} \label{R}
\sum_1^\infty a_n^{-2}= R^2 <  \infty
\end{equation}
and let $x_n \in H$
be vectors with $\|x_n\|=a_n$. We may assume that $H$ is a
complex Hilbert space, since otherwise we may embed $H$ into
its complexification. Define planks
$$
P_n=\{h \in H : |\langle h \, , \, x_n \rangle| \le \frac{1}{2}\}.
$$
The width of $P_n$ equals $a_n^{-1}$. Using (\ref{R}) and
the Theorem \ref{Ball} we deduce that the planks $P_n$
cannot cover the whole space $H$ (they even cannot cover
a ball of radius $R + \eps$). So there is an element
$h \in H$ for which all the inequalities
$$
|\langle h \, , \, x_n \rangle| > \frac{1}{2}
$$
hold true at the same time. This $h$ separates our sequence
$x_n$ from 0.
\end{proof}
\section{A comment to the Ball's Theorem \ref{Ball}}
Why estimates for coverings by real and complex planks
differ so strongly?
A possible explanation looks as follows: a complex plank
of width $r$ looks not like a slice between two hyperplanes,
but like an orthogonal cylinder having as a base a circle
of the radius $r/2$. If one tries to cover a circle of radius
$R$ by circles of radiuses $r_n$, then calculating corresponding
areas one can easily see, that $\sum r_n^2 \ge R^2$.
This explanation leads to the following hypothesis:
let $k \in \N$ be a fixed number, and let $C_n$ be orthogonal
cylinders in a real Hilbert space $H$, having as their bases
$k$-dimensional balls of
radiuses $r_n$ respectively and $\bigcup_{n \in \N} C_n = H$.
Then $\sum_{n=1}^\infty r_n^k = \infty$. With the help of the
Theorem \ref{main}
we can show that this natural hypothesis fails already for $k=3$.
In fact, assume that the abovementioned hypothesis is true for $k=3$.
Consider a sequence $x_n \in H$, $\|x_n\|=a_n$, having 0 as a weak
cluster point, but with $\sum_1^\infty a_n^{-3} <  \infty$ (such a
sequence exists due to the Theorem \ref{main}). Consider an auxiliary
Hilbert space $H_1 = H \oplus H \oplus H$ - the orthogonal direct sum
of 3 copies of the original space $H$. Introduce cylinders $C_n$ in
$H_1$ as follows:
$$
C_n=\{h=(h_1,h_2,h_3) \in H_1 :
\sum_{j=1}^3 | \langle h_j \, , \, x_n \rangle |^2 \le 1\}.
$$
$C_n$ are orthogonal
cylinders in $H_1$, having as their bases
$3$-dimensional balls of radiuses $a_n^{-1}$.
According to our hypothesis, these cylinders do not cover the whole
space $H_1$, i.e. there is an $g=(g_1,g_2,g_3) \in H_1$ which does not
belong to any of these cylinders. But this means, that the weak
neighborhood of 0
$$
W =\{x \in H :
\sum_{j=1}^3 | \langle g_j \, , \, x \rangle |^2 \le 1\}
$$
separates all the $x_n$ from 0. Contradiction.
\section{Generalization to Banach spaces:
the role of finite representability and cotype}
Let $X$ be a Banach space, $2 \le p \le \infty$. The space $l_p$
is finitely representable in $X$ if for every $\eps > 0$
and for every $n \in \N$ there are elements $e_1, e_2, \ldots e_n \in S_X$
such that
$$
(1-\eps) \left ( \sum_{j=1}^n |b_j|^p \right )^\frac{1}{p}
\le \| \sum_{j=1}^n b_j e_j\|
\le (1+\eps) \left ( \sum_{j=1}^n |b_j|^p \right )^\frac{1}{p}
$$
for all selections of coefficients $b_j$.
\begin{thm} \label{fin-p}
Let $X$ be a Banach space, $2 \le p \le \infty$ and
let $1 \le p' \le 2$ be dual to $p$ exponent. Let moreover
$l_p$ be finitely representable in $X$ and
$a_n>0$ satisfy condition $\sum_1^\infty a_n^{-p'}=\infty$.
Then there is a sequence $x_n \in X$ with
$\|x_n\|=a_n$, having 0 as a weak cluster point. Moreover
there is a matrix $P=\{p_{n,m}\}$ of weights
such that $w(P,p')\lim x_n =0$.
\end{thm}
\begin{proof} Fix an $\eps <  1/8$.
First choose $0=n_1 <  n_2 <  n_3 <  \ldots$ to satisfy condition
$$
\lim_{k \to \infty}\sum_{j=n_k+1}^{n_{k+1}} a_j^{-p'}=\infty.
$$
Define $P=\{p_{k,m}\}$ as follows:
$$p_{k,m}=\frac{a_m^{-p'}}{\sum_{j=n_k+1}^{n_{k+1}} a_j^{-p'}}$$
when $n_k<  m \le n_{k+1}$ and $p_{k,m}=0$ otherwise.
Using step-by-step finite reperesentability of $l_p$ in $X$
select sequence $e_n \in S_X$, for which
$$
(1-\eps) \left ( \sum_{j=n_k+1}^{n_{k+1}} |b_j|^p \right )^\frac{1}{p}
\le \| \sum_{j=n_k+1}^{n_{k+1}} b_j e_j\|
\le (1+\eps) \left ( \sum_{j=n_k+1}^{n_{k+1}} |b_j|^p \right )^\frac{1}{p}
$$
for all $k$ and $b_j$.
Now define $x_n = a_n e_n$. Then for every $f \in X^*$ we have
$$
\sum_{m=1}^\infty p_{k,m} \left | f(x_m) \right |^{p'}=
\frac{\sum_{j=n_k+1}^{n_{k+1}}
\left | f(e_m) \right |^{p'}}{\sum_{j=n_k+1}^{n_{k+1}} a_j^{-p'}} \le
\frac{2\|f\|^{p'}}{\sum_{j=n_k+1}^{n_{k+1}} a_j^{-p'}} \to 0,
\ as \ n \to \infty
$$
which means $w(P,p')$-convergence of $x_n$ to 0.
By the Proposition \ref{wp-conv} this means that 0 is a weak
cluster point
of $\{x_n\}$.
\end{proof}
Together with Dvoretzky's theorem ($l_2$ is finitely representable in
every infinite-dimensional Banach space) this gives us the following:
\begin{cor}
For every infinite-dimensional Banach space $X$
and every selection of $a_n>0$ satisfying condition
$\sum_1^\infty a_n^{-2}=\infty$
there is a sequence $x_n \in X$ with
$\|x_n\|=a_n$, having 0 as a weak cluster point.
\end{cor}
Recall, that a Banach space $X$ has M-cotype $p <  \infty$ if there
is
a constant $C > 0$ such that for every finite collection of vectors
$\{x_k\}_{k=1}^n \subset X$ there are coefficients
$\gamma_k = \pm 1$ for which
$$
\| \sum_{k=1}^n \gamma_k x_x \| \ge
C \left( \sum_{k=1}^n \|x_k \|^p \right)^{1/p}.
$$
Due to Maurey - Pisier theorem a space $X$ has an M-cotype
if and only if $l_\infty$ is not finitely
representable in $X$. A small survey of facts concerning M-cotype
can be found in \cite{kadkad}.

The Theorem \ref{fin-p} together with another Ball's result \cite{ball1}
(generalization of the Bang's real plank theorem to arbitrary Banach
spaces) gives us an analog of Theorem \ref{main} for spaces without
cotype (in particular for $c_0$). In this case the characterization
does not involve square exponentials.
\begin{cor}
Let $X$ be a Banach spaces in which $l_\infty$ is finitely
representable. Then the following properties for a sequence
$a_n>0$ are equivalent:
\begin{enumerate}
\item \label{dva2}
There is a sequence $x_n \in X$ with
$\|x_n\|=a_n$, and there is a sequence $P=\{p_{n,m}\}$ of weights
such that $w(P,1)\lim x_n =0$.
\item There is a sequence $x_n \in X$ with
$\|x_n\|=a_n$, having 0 as a weak cluster point;
\item
$\sum_1^\infty a_n^{-1}=\infty$.
\end{enumerate}
\end{cor}
The next theorem gives us a relationship between weak
weighted limit and cotype.
\begin{thm} \label{cotype}
Let $X$ be a Banach space of M-cotype $2 \le p \le \infty$ and
let $1 \le p' \le 2$ be dual to $p$ exponential. Let $x_n \in X$
be a sequence
such that $w(P,p')\lim x_n =0$
for some matrix $P=\{p_{n,m}\}$ of weights. Then
$\sum_1^\infty \|x_n\|^{-p'}=\infty$.
\end{thm}
\begin{proof}
Since $\lim_{n \to \infty} p_{n,m} =0$, using small perturbation
argument we may assume that there are 
$m_1(n) <  m_2(n)$, $m_1(n) \to \infty$,
such that $p_{n,m} =0$ for $m$ outside the interval $(m_1(n) , m_2(n))$.
By the closed graph theorem, applied to the operator
$T: X^* \to \left ( \sum_{n=1}^\infty l_1\right )_\infty $,
$$
Tx^* = \bigl ( (p_{1,1}|x^*(x_1)|,p_{1,2}|x^*(x_2)| \ldots); \
(p_{2,1}x^*(x_1),p_{2,2}x^*(x_2) \ldots); \ \ldots \bigr ),
$$
there is a constant $C>0$ such that
$$
\sum_{m=1}^\infty p_{n,m}|x^*(x_m)| \le C \| x^*\|
$$
for all $n \in \N$ and all $x^* \in X^*$. So
$$
C \ge \sup_{x^* \in S_{X^*}} \sum_{m=1}^\infty p_{n,m}|x^*(x_m)|=
\sup_{x^* \in S_{X^*}} \sum_{m=m_1(n)}^{m_2(n)} p_{n,m}|x^*(x_m)|=
$$
$$
\sup_{x^* \in S_{X^*}} \sup_{\gamma_m= \pm 1}
\sum_{m=m_1(n)}^{m_2(n)} p_{n,m}\gamma_mx^*(x_m) =
\sup_{\gamma_m = \pm 1}
\| \sum_{m=m_1(n)}^{m_2(n)} p_{n,m}x_m \| \ge
$$
$$
C_1 \left ( \sum_{m=m_1(n)}^{m_2(n)} (p_{n,m}\|x_m\|)^p \right )^{1/p},
$$
where $C_1$ is the constant from the definition of M-cotype.
Applying H\"older inequality to $\sum_{m=m_1(n)}^{m_2(n)} p_{n,m}=1$
we deduce
$$
1 \le \left ( \sum_{m=m_1(n)}^{m_2(n)} (p_{n,m}\|x_m\|)^p \right )^{1/p}
\left ( \sum_{m=m_1(n)}^{m_2(n)} \|x_m\|^{-p'} \right )^{1/p'}
\le \frac{C}{C_1} \left ( \sum_{m=m_1(n)}^{m_2(n)} \|x_m\|^{-p'} \right)^{1/p'}
$$
which means that $\sum_1^\infty \|x_n\|^{-p'}=\infty$.
\end{proof}


\begin{thebibliography}{10}
\bibitem{ball} Keith Ball, {\it Convex Geometry and Functional Analysis},
in W.B.Johnson and J.Lindenstrauss (editors) Handbook of the geometry
of Banach spaces, vol. 1 (2001), 161 - 194.
\bibitem{ball1} Keith Ball, {\it The plank problem for symmetric bodies},
Invent. Math. {\bf 104} (1991), 535 - 543.
\bibitem{comp} Keith Ball, {\it The complex plank problem},
Bull. London Math. Soc. {\bf 33} (2001) no. 4, 433 - 442.
\bibitem{bang} T. Bang, {\it A solution of the "Plank problem"}, Proc.
Amer. Math. Soc. {\bf 2} (1951), 990 - 993
\bibitem{conkad} J.Connor, M.Ganichev, V.Kadets {\it A characterization of
Banach spaces with separable duals via weak statistical convergence},
J. Math. Anal. Appl. {\bf 244} (2000), no 1, 251 - 261.
\bibitem{kadkad} M.Kadets, V.Kadets
{\it Series in Banach spaces. Conditional and unconditional convergence.}
Translated from the Russian by Andrei Iacob.
Operator Theory: Advances and Applications, 94. Birkhauser Verlag,
Basel 1997.
\end{thebibliography}
\end{document}